\documentclass[letter, 10pt, conference]{cssconf}
\IEEEoverridecommandlockouts

\usepackage{amsfonts,latexsym,graphicx,epsfig,amssymb,amsmath}
\usepackage{color}
\usepackage{cite}
\usepackage{graphicx}
\usepackage[bookmarks=false,bookmarksopen=false,colorlinks,    citecolor={blue}, filecolor={blue}, linkcolor={blue}, urlcolor={blue}]{hyperref}

\begin{document}

\allowdisplaybreaks

    \newtheorem{thm}{Theorem}[section]

    \newtheorem{cor}{Corollary}[section]

    \newtheorem{lem}{Lemma}[section]

    \newtheorem{prop}{Proposition}[section]

    \newtheorem{fact}{Fact}[section]

    \newtheorem{claim}{Claim}[section]

    \newtheorem{prb}{Problem}[section]
    \newtheorem{rem}{Remark}[section]
    \newtheorem{definition}{Definition}[section]

    \newtheorem{assume}{Assumption}[section]

    \newtheorem{eg}{Example}[section]

    \newtheorem{ceg}{Counter Example}[section]

    \newenvironment{pf}{\noindent \textbf{Proof.}}{\hfill{$\blacksquare$}}

    \newenvironment{note}{Note:}

\newcommand{\Jcal}{\mathcal{J}}
\newcommand{\Rcal}{\mathcal{R}}

\newcommand{\etab}{\boldsymbol{\eta}}
\newcommand{\varsigmab}{\boldsymbol{\varsigma}}
\newcommand{\phib}{\boldsymbol{\varphi}}
\newcommand{\alphab}{\boldsymbol{\alpha}}
\newcommand{\betab}{\boldsymbol{\beta}}
\newcommand{\gammab}{\boldsymbol{\gamma}}
\newcommand{\psib}{\boldsymbol{\psi}}
\newcommand{\mub}{\boldsymbol{\mu}}
\newcommand{\nub}{\boldsymbol{\nu}}
\newcommand{\xib}{\boldsymbol{\xi}}
\newcommand{\omegab}{\boldsymbol{\omega}}
\newcommand{\pib}{\boldsymbol{\pi}}
\newcommand{\taub}{\boldsymbol{\tau}}
\newcommand{\lambdab}{\boldsymbol{\lambda}}
\newcommand{\deltab}{\boldsymbol{\delta}}

\newcommand{\Xb}{\mathbf{X}}
\newcommand{\Yb}{\mathbf{Y}}
\newcommand{\Ab}{\mathbf{A}}
\newcommand{\Cb}{\mathbf{C}}
\newcommand{\db}{\mathbf{d}}
\newcommand{\zb}{\mathbf{z}}
\newcommand{\qb}{\mathbf{q}}
\newcommand{\Db}{\mathbf{D}}
\newcommand{\Rb}{\mathbf{R}}
\newcommand{\Mb}{\mathbf{M}}
\newcommand{\Bb}{\mathbf{B}}
\newcommand{\Wb}{\mathbf{W}}
\newcommand{\Jb}{\mathbf{J}}
\newcommand{\Ib}{\mathbf{I}}
\newcommand{\Vb}{\mathbf{V}}
\newcommand{\Ub}{\mathbf{U}}
\newcommand{\gb}{\mathbf{g}}
\newcommand{\fb}{\mathbf{f}}
\newcommand{\Hb}{\mathbf{H}}

\newcommand{\Omegab}{\mathbf{\Omega}}
\newcommand{\Thetab}{\mathbf{\Theta}}
\newcommand{\Pib}{\mathbf{\Pi}}
\newcommand{\Lambdab}{\mathbf{\Lambda}}
\newcommand{\Deltab}{\mathbf{\Delta}}

\newcommand{\Trm}{\textrm{T}}
\newcommand{\Psf}{\textsf{P}}
\newcommand{\Qsf}{\textsf{Q}}
\newcommand{\Tsf}{\textsf{T}}
\newcommand{\TPsf}{\textsf{TP}}
\newcommand{\TsPsf}{\textsf{T}^*\textsf{P}}
\newcommand{\TQsf}{\textsf{TQ}}
\newcommand{\TsQsf}{\textsf{T}^*\textsf{Q}}
\newcommand{\Gsf}{\textsf{G}}
\newcommand{\TsGsf}{\textsf{T}^*\textsf{G}}
\newcommand{\SO}{\textsf{SO}(3)}
\newcommand{\TSO}{\textsf{TSO}(3)}
\newcommand{\SOtwo}{\textsf{SO}(2)}
\newcommand{\SEthree}{\textsf{SE}(3)}
\newcommand{\SU}{\textsf{SU}(N)}

\newcommand{\gfk}{\mathfrak{g}}
\newcommand{\so}{\mathfrak{so}}
\newcommand{\se}{\mathfrak{se}}

\newcommand{\drm}{{\rm d}}
\newcommand{\Drm}{{\rm D}}
\newcommand{\Ad}{{\rm Ad}}
\newcommand{\ad}{{\rm ad}}

\newcommand{\be}{\begin{eqnarray}}
\newcommand{\ee}{\end{eqnarray}}

\newcommand{\nn}{\nonumber}

\newcommand{\beas}{\begin{eqnarray*}}
\newcommand{\eeas}{\end{eqnarray*}}
\newcommand{\la}{\label}
\newcommand{\T}{^{\mbox{\small T}}}

\abovedisplayskip=0.1ex 
\belowdisplayskip=0.1ex 
\abovedisplayshortskip=0.1ex 
\belowdisplayshortskip=0.1ex 

\title{\LARGE A Discrete Variational Integrator for Optimal Control Problems on
$\SO$}

\author{Islam I. Hussein\thanks{Islam Hussein is Assistant Professor of Mechanical Engineering at Worcester Polytechnic Institute, \href{mailto:ihussein@wpi.edu}
{ihussein@uiuc.edu}.} \and Melvin Leok\thanks{Melvin Leok is Assistant Professor of Mathematics at Purdue University, West Lafayette, \href{mailto:mleok@math.purdue.edu}
{mleok@umich.edu}.} \and Amit K. Sanyal\thanks{Amit Sanyal is
Post-Doctoral Research Associate at Arizona State University, Tempe,
\href{mailto:sanyal@asu.edu} {sanyal@asu.edu}.} \and Anthony M.
Bloch\thanks{Anthony Bloch is Professor of Mathematics at the
University of Michigan, Ann Arbor, \href{mailto:abloch@umich.edu}
{abloch@umich.edu}.}}

\maketitle

\begin{abstract}

In this paper we study a discrete variational optimal control
problem for the rigid body. The cost to be minimized is the external
torque applied to move the rigid body from an initial condition to a
pre-specified terminal condition. Instead of discretizing the
equations of motion, we use the discrete equations obtained from the
discrete Lagrange--d'Alembert principle, a process that better
approximates the equations of motion. Within the discrete-time
setting, these two approaches are not equivalent in general. The
kinematics are discretized using a natural Lie-algebraic formulation
that guarantees that the flow remains on the Lie group $\SO$ and its
algebra $\so(3)$. We use Lagrange's method for constrained problems
in the calculus of variations to derive the discrete-time necessary
conditions. We give a numerical example for a three-dimensional
rigid body maneuver.

\end{abstract}

\section{Introduction}\label{sec:intro}

This paper deals with a structure-preserving computational approach
to the optimal control problem of minimizing the control effort
necessary to perform an attitude transfer from an initial state to a
prescribed final state, in the absence of a potential field. The
configuration of the rigid body is given by the rotation matrix from
the body frame to the spatial frame, which is an element of the
group of orientation-preserving isometries in $\mathbb{R}^3$. The
state of the rigid body is described by the rotation matrix and its
angular velocity.

To motivate the computational approach we adopt in the discrete-time
case, we first revisit the variational continuous-time optimal
control problem. The continuous-time extremal solutions to this
optimal control problem have certain special features, since they
arise from variational principles. General numerical integration
methods, including the popular Runge-Kutta schemes, typically
preserve neither first integrals nor the characteristics of the
configuration space. Geometric integrators are the class of
numerical integration schemes that preserve such properties, and a
good survey can be found in \cite{Hairer:02}.
Techniques particular to Hamiltonian systems are also discussed in
\cite{Leimkuhler:04} and \cite{SanzSerna:94}.

Our approach to discretizing the optimal control problem is in
contrast to traditional techniques such as collocation, wherein the
continuous equations of motion are imposed as constraints at a set
of collocation points. In our approach, modeled after
\cite{Junge:05}, the discrete equations of motion are derived from a
discrete variational principle, and this induces constraints on the
configuration at each discrete time step.

This approach yields discrete dynamics that are more faithful to the
continuous equations of motion, and consequently yields more
accurate numerical solutions to the optimal control problem. This feature is extremely important in computing
accurate (sub)optimal trajectories for long-term spacecraft attitude
maneuvers. For example, in \cite{hussein:03c}, the authors propose
an imaging spacecraft formation design that requires a continuous
attitude maneuver over a period of 77 days in a low Earth orbit.
Hence, the attitude maneuver has to be very accurate to meet tight
imaging constraints over long time ranges. The proposed variational
scheme can also be easily extended to other types of Lie groups. For
example, in long range inter-planetary orbit transfers (see, for
example, \cite{NewHorizons:06}), one is interested in computing
optimal or suboptimal trajectories on the group of rigid body
motions $\SEthree$ with a high degree of accuracy. Similar
requirements also apply to the control of quantum systems. For
example, efficient construction of quantum gates is a problem on the
unitary Lie group $\SU$. This is an optimal control problem, where
one wishes to steer the identity operator to the  desired unitary
operator (see, for example, \cite{PalaoKosloff:02} and
\cite{KhanejaGlaserBrockett:02}).

Moreover, an important feature of the way we discretize the optimal
control problem is that it is $\SO$-equivariant. This is desirable, since it
ensures that our numerical results are independent of the choice of coordinates
and coordinate frames. This is in contrast to methods based on
coordinatizing the rotation group using quaternions, (modified)
Rodrigues parameters, and Euler angles, as given in the survey
\cite{ScriTho94}. Even if the optimal cost function is
$\SO$-invariant, as in \cite{ScJuRo96}, the use of generalized
coordinates imposes constraints on the attitude kinematics.

For the purpose of numerical simulation, the corresponding discrete
optimal control problem is posed on the discrete state space as a
two stage discrete variational problem. In the first step, we derive
the discrete dynamics for the rigid body in the context of discrete
variational mechanics \cite{marsden:01}. This is achieved by
considering the discrete Lagrange--d'Alembert variational principle
\cite{Kane:00} in combination with essential ideas from Lie group
methods \cite{Iserles:00}, which yields a Lie group variational
integrator \cite{Leok:04}. This integrator explicitly preserves the
Lie group structure of the configuration space, and is similar to
the integrators introduced in \cite{Lee:05} for a rigid body in an
external field, and in \cite{Lee:05b} for full body dynamics. These
discrete equations are then imposed as constraints to be satisfied
by the extremal solutions to the discrete optimal control problem,
and we obtain the discrete extremal solutions in terms of the given
terminal states.

The paper is organized as follows. As motivation, in Section
\ref{sec:continuous_time}, we study the optimal control problem in
continuous-time. In Section \ref{sec:discrete_time}, we study the
discrete-time case. In particular, in Section
\ref{subsec:formulation_discrete_example} we state the optimal
control problem and describe our approach. In Section
\ref{sec:discrete_lag_RB}, we derive the discrete-time equations of
motion for the rigid body starting with the discrete Lagrange--d'Alembert principle. These equations are used in Section
\ref{subsec:Discrete_Variational} for the optimal control problem.
In Section \ref{sec:numerics}, we describe an algorithm for solving
the general nonlinear, implicit necessary conditions for $\SO$ and
give numerical examples for rest-to-rest and slew-up spacecraft maneuvers.

\section{Continuous-Time Results}\label{sec:continuous_time}

\subsection{Problem
Formulation}\label{subsec:formulation_continuous_example}

In this paper, the
natural pairing between $\so^*(3)$ and $\so(3)$ is denoted by
$\left<\cdot,\cdot\right>$. Let $\ll\cdot,\cdot\gg$ and $\ll\cdot,\cdot\gg_*$ denote
the standard (induced by the Killing form) inner product on $\so(3)$ and $\so^*(3)$,
respectively. The inner product $\ll\cdot,\cdot\gg_*$ is naturally
induced from the standard norm $\ll\xib,\omegab\gg=-\frac{1}{2}{\rm Tr}(\xib^T\omegab)$, for all $\xib,\omegab\in\so(3)$,
through\be\label{eq:inner_products}\ll\etab,\phib\gg_*&=&\left<\etab,\phib^{\sharp}\right>=\left<\etab,\omegab\right>=\left<\xib^{\flat},\omegab\right>\nn\\&=&\ll\xib,\omegab\gg,\ee where $\phib=\omegab^{\flat}\in\so^*(3)$ and
$\etab=\xib^{\flat}\in\so^*(3)$, with $\xib,\omegab\in\so(3)$ and $\flat$ and $\sharp$ are the musical isomorphisms (see \S2.5 of \cite{bloch:03}) with respect to the standard metric $\ll\cdot,\cdot\gg$. On $\so(3)$, these isomorphisms correspond to the transpose operation. That is, we have $\phib=\omegab^{\Trm}$ and
$\etab=\xib^{\Trm}$.

Let $\Jb:\so(3)\rightarrow\so^*(3)$ be the positive definite inertia operator. It can be shown that
\be\label{eq:inertia_property}\left<\Jb(\xib),\omegab\right>=\left<\Jb(\omegab),\xib\right>.\ee On $\so(3)$, $\Jb$ is given by $\Jb(\xib)=J\xib+\xib J$, where
$J$ is a positive definite symmetric matrix (see, for example, \cite{bloch:03,hussein:CDC05}). Moreover, we also have
$\Jb(\etab^{\sharp})^{\sharp}=(J\etab^{\Trm}+\etab^{\Trm}J)^{\Trm}=\Jb(\etab)$,
which is an abuse of notation since $\etab\in\so^*(3)$. For the sake of generality and mathematical accuracy, we will use the general definitions, though it is helpful to keep the above identifications
for $\so(3)$ in mind.

In this section we review some continuous-time optimal control
results using a simple optimal control example on $\SO$. Here, we minimize
the norm squared of the control
torque $\taub\in\so^*(3)$ applied to rotate a rigid body subject to the Lagrange--d'Alembert principle for the rigid
body\footnote{This is equivalent to constraining the problem to
satisfy the rigid body equations of motion given by
equations (\ref{eq:SO3_eoms_example}). However, for the sake
of generality that will be appreciated in the discrete-time problem,
we choose to treat the Lagrange--d'Alembert principle as the
constraint as opposed to the rigid body equations of motion. Both
are equivalent in the continuous-time case but are generally not equivalent in the discrete-time case.}
whose configuration is given by $\Rb\in\SO$ and body angular
velocity is given by $\Omegab\in\so(3)$. We require that the system
evolve from an initial state $(\Rb_0,\Omegab_0)$ to a final state
$(\Rb_T,\Omegab_T)$ at a fixed terminal time $T$. Hence, we have the
following minimum control effort optimal control problem.

\begin{prb}\label{pb:SO3_example0} Minimize
\be\label{eq:SO3_cost_LD_example0}
  \Jcal &=&
  \frac{1}{2}\int_0^T\ll\taub,\taub\gg_*\drm t
\ee subject to
\begin{enumerate}
    \item satisfying \emph{Lagrange--d'Alembert principle}:
     \be\label{eq:SO3_LD_example0}
      \deltab\int_0^T\frac{1}{2}\left<\Jb\left(\Omegab\right),\Omegab\right>\drm t+\int_0^T\left<\taub,\Wb\right>\drm t=0,
    \ee
    subject to $\dot{\Rb}=\Rb\Omegab$, where
    $\Wb$ is the variation vector field to be defined below,
    \item and the \emph{boundary conditions} \be\label{eq:SO3_bcs_example0}
    \Rb(0)&=&\Rb_0, ~\Omegab(0)=\Omegab_0,\nn\\ \Rb(T)&=&\Rb_T,
    ~\Omegab(T)=\Omegab_T.\ee
\end{enumerate}
\end{prb}

We now show that this is equivalent to
the following problem formulation,
where the rigid body equations of motion replace the Lagrange--d'Alembert principle.

\begin{prb}\label{pb:SO3_example} Minimize
\be\label{eq:SO3_cost_example}
  \Jcal &=&
  \frac{1}{2}\int_0^T\ll\taub,\taub\gg_*\drm t
\ee subject to
\begin{enumerate}
    \item the \emph{dynamics} \be\label{eq:SO3_eoms_example}
      \dot{\Rb} &=& \Rb\Omegab\\
      \dot{\Mb} &=&\ad_{\Omegab}^*\Mb+\taub=
      \left[\Mb,\Omegab\right]+\taub,\nn\ee
      where $\Mb=\Jb(\Omegab)\in\so^*(3)$ is the momentum,
    \item and the \emph{boundary conditions} \be\label{eq:SO3_bcs_example}
    \Rb(0)&=&\Rb_0, ~\Omegab(0)=\Omegab_0,\nn\\ \Rb(T)&=&\Rb_T, ~\Omegab(T)=\Omegab_T.
\ee
\end{enumerate}
\end{prb}

In the above, $\ad^*$ is the dual of the adjoint representation,
$\ad$, of $\so(3)$ and is given by
$\ad_{\xib}^*\etab=-[\xib,\etab]\in\so^*(3)$, for all $\xib\in\so(3)$ and
$\etab\in\so^*(3)$. Recall that the bracket is defined by
$\left[\xib,\omegab\right]=\xib\omegab-\omegab\xib$.

\subsection{The Lagrange--d'Alembert Principle and the Rigid
Body Equations of Motion}\label{subsec:LD_continuous_example}

In this section we derive the forced rigid body equations of motion
(equations (\ref{eq:SO3_eoms_example})) from the Lagrange--d'Alembert
principle. In dealing with the kinematic constraint, 
$\dot{\Rb}=\Rb\Omegab$, we may either append it to the Lagrangian using the method of Lagrange multipliers, or we can directly compute the constrained variations (see \S13.5 of \cite{marsden:99}).  Here, we take the direct approach as it yields a more concise derivation.

First, we take
variations of the kinematic condition $\Omegab=\Rb^{-1}\Rb$ to
obtain $\deltab\Omegab =
  -\Rb^{-1}\left(\deltab\Rb\right)\Rb^{-1}\dot{\Rb}+\Rb^{-1}\deltab\dot{\Rb}$. As defined previously, we have $\Wb=\Rb^{-1}\deltab\Rb$ and,
therefore, $\dot{\Wb}=-\Rb^{-1}\dot{\Rb}\Rb^{-1}\deltab\Rb+\Rb^{-1}\deltab\dot{\Rb}=-\Omegab\Wb+\Rb^{-1}\deltab\dot{\Rb}$,
since $\deltab\dot{\Rb}=\frac{\drm}{\drm t}\deltab\Rb$ (see \cite{milnor:63}, p. 52). Hence, we have
\be\label{eq:direct_variation_cont_example}
  \deltab\Omegab =
  -\Wb\Omegab+\Omegab\Wb+\dot{\Wb}
  =\ad_{\Omega}\Wb+\dot{\Wb}.
\ee Taking variations of the Lagrange--d'Alembert principle we obtain \be
    \int_0^T\left<\Jb\left(\Omegab\right),\deltab\Omegab\right>+\left<\taub,\Wb\right>\drm t=0.\nn
\ee Using the variation in equation
(\ref{eq:direct_variation_cont_example}) and integrating by parts,
we obtain \be
    0=\int_0^T\left<-\dot{\Mb}+\ad_{\Omega}^*\Mb+\taub,\Wb\right>\drm
    t+\left[\left<\Jb\left(\Omegab\right),\Wb(t)\right>\right]_0^T,\nn
\ee where we have used the
property\be\label{eq:ad_property}\left<\etab,\ad_{\omegab}\xib\right>=\left<\ad_{\omegab}^*\etab,\xib\right>,
~\etab\in\so^*(3), ~\omegab,\xib\in\so(3).\ee This gives the desired
result, with $\Mb=\Jb\left(\Omegab\right)$.

\subsection{Continuous-Time Variational Optimal Control Problem}\label{subsec:Continuous_Variational}

A direct variational approach is used here to derive the necessary conditions for the optimal control Problem
(\ref{pb:SO3_example}).

\emph{\textbf{A Second Order Direct Approach.}} ``Second order" is
used here to reflect the fact that we now study variations of second
order dynamical equations as opposed to the kinematic direct
approach studied in Section \ref{subsec:LD_continuous_example}. We
now give the resulting necessary conditions using a direct approach
as in \cite{marsden:99}. We already computed the variations of $\Rb$
and $\Omegab$. These were as follows: $\deltab\Rb=\Rb\Wb$ and
$\deltab\Omegab=\ad_{\Omegab}\Wb+\dot{\Wb}$. We now compute the
variation of $\dot{\Mb}$ with the goal of obtaining the proper
variations for $\taub$: \be
  \deltab\dot{\Mb} &=&
  \Jb\left(\deltab\dot{\Omegab}\right)=\Jb\left(\frac{\drm}{\drm
  t}\deltab\Omegab+\Rcal\left(\Wb,\Omegab\right)\Omegab\right),\nn
\ee where $\Rcal$ is the curvature tensor on SO(3). The curvature
tensor $\Rcal$ arises due to the identity (see \cite{milnor:63},
page 52) \be\label{eq:curvature}
    \frac{\partial}{\partial\epsilon}\frac{\partial}{\partial t}\mathbf{Y}-
    \frac{\partial}{\partial t}\frac{\partial}{\partial
    \epsilon}\mathbf{Y}=\Rcal(\Wb,\mathbf{Y})\Omegab,\nn
\ee where $\mathbf{Y}\in\TSO$ is any vector field along the curve
$\Rb(t)\in\SO$. Taking variations of
$\dot{\Mb}=\ad_{\Omegab}^*\Mb+\taub$ we obtain $\deltab\dot{\Mb}
  =\ad_{\deltab\Omegab}^*\Mb+\ad_{\Omegab}^*\deltab\Mb+
  \deltab\taub$. We now have the desired variation in $\taub$:
\be\la{eq:tau_variations}
  \deltab\taub &=&
  \Jb\left(\Rcal\left(\Wb,\Omegab\right)\Omegab\right)+\frac{\drm }{\drm
  t}\Jb\left(\deltab\Omegab\right)-\ad_{\deltab\Omegab}^*\Mb\nn\\&&-\ad_{\Omegab}^*\deltab\Mb.
\ee
Take variations of the cost functional
(\ref{eq:SO3_cost_example}) to obtain: \be
  &&\deltab\Jcal =   \int_0^T\big<\Jb(\ddot{\varsigma})-\ad_{\Omegab}^*\left(\Jb(\dot{\varsigma})\right)+\dot{\etab}-\frac{\drm}{\drm
  t}\left(\ad_{\varsigmab}^*\Mb\right)\nn\\
  &&+\left[\Rcal\left(\Jb(\varsigmab)^{\sharp},\Omegab\right)\Omegab\right]^{\flat}+\ad_{\Omegab}^*\ad_{\varsigmab}^*\Mb
  -\ad_{\Omegab}^*\etab,\Wb\big>\drm
  t,\nn
\ee where $\varsigmab=\taub^{\sharp}\in\so(3)$ and
$\etab=\Jb\left(\ad_{\Omegab}\varsigmab\right)\in\so^*(3)$. Here, we used integration by parts
and the boundary conditions (\ref{eq:SO3_eoms_example}), equations
(\ref{eq:direct_variation_cont_example}) and
(\ref{eq:tau_variations}), and the identities
(\ref{eq:inner_products}), (\ref{eq:inertia_property}) and
(\ref{eq:ad_property}). Hence, we have the following theorem.

\begin{thm}\label{thm:SO3_example2} The necessary optimality conditions for the problem of minimizing
(\ref{eq:SO3_cost_example}) subject to the dynamics
(\ref{eq:SO3_eoms_example}) and the boundary conditions
(\ref{eq:SO3_bcs_example}) are given by the single fourth
order\footnote{Second order in $\taub$ and fourth order in $\Rb$.}
differential equation \be
  0 &=& \Jb(\ddot{\varsigma})-\ad_{\Omegab}^*\left(\Jb(\dot{\varsigma})\right)+\dot{\etab}-\frac{\drm}{\drm
  t}\left(\ad_{\varsigmab}^*\Mb\right)\nn\\
  &&+\left[\Rcal\left(\left(\Jb(\varsigmab)\right)^{\sharp},\Omegab\right)\Omegab\right]^{\flat}+\ad_{\Omegab}^*\left(\ad_{\varsigmab}^*\Mb\right)
  -\ad_{\Omegab}^*\etab,\nn
\ee as well as the equations (\ref{eq:SO3_eoms_example}) and the
boundary conditions (\ref{eq:SO3_bcs_example}), where $\varsigmab$
and $\etab$ are as defined above.
\end{thm}

To obtain above result we used the initial conditions
(\ref{eq:SO3_bcs_example}), and the fact that the vector
fields $\Omegab$ and $\Wb$ are left-invariant vector fields. The
curvature tensor is evaluated at a point $\Rb(t)\ne \Ib$. That is,
we get $\Rcal_{\Rb_{\epsilon}}\left(\frac{\partial
\Rb_{\epsilon}}{\partial\epsilon},\frac{\partial
\Rb_{\epsilon}}{\partial t}\right)\Omegab$. Evaluating this at
$\epsilon=0$ we get:
$\Rcal_{\Rb}\left(\Rb\Wb,\Rb\Omegab\right)\Omegab$. Since $\Rb\Wb$
and $\Rb\Omegab$ are left-invariant vector fields at the group
element $\Rb(t)$, by the identification $\Tsf_{\Rb}\SO\simeq\so(3)$,
we have
$\Rcal_{\Rb}\left(\Rb\Wb,\Rb\Omegab\right)\Omegab=\Rcal\left(\Wb,\Omegab\right)\Omegab$,
which is the curvature tensor evaluated at the identity element. For a compact semi-simple Lie group $\Gsf$ with Lie algebra
$\mathfrak{g}$, the curvature tensor, with respect to a bi-invariant
metric, is (see \cite{milnor:63}):
\be\label{eq:curvature_compact}
\Rcal\left(\mathbf{X},\mathbf{Y}\right)\mathbf{Z}=\frac{1}{4}\ad_{\ad_{\Xb}\Yb}\mathbf{Z},
\ee for all $\mathbf{X},\mathbf{Y},\mathbf{Z}\in\mathfrak{g}$.

Using a Lagrange multiplier approach, we obtain instead the following theorem.

\begin{thm}\label{thm:SO3_example} The necessary optimality conditions
for the problem of minimizing (\ref{eq:SO3_cost_example}) subject to
the dynamics (\ref{eq:SO3_eoms_example}) and the boundary conditions
(\ref{eq:SO3_bcs_example}) are given by
\be\label{eq:SO3_theorem_example}
  \taub &=& \Lambdab_2\nn \\
\dot{\Lambdab}_1&=&\left[\Rcal\left(\Jb\left(\Lambdab_2\right)^{\sharp},\Omegab\right)\Omegab\right]^{\flat}+\ad_{\Omegab}^*\Lambdab_1\\
\dot{\Lambdab}_2&=&-\Jb^{-1}\left(\Lambdab_1\right)-\ad_{\Omegab}\Lambdab_2+\Jb^{-1}\left(\ad_{\Lambdab_2}^*\Mb\right)\nn
\ee and the equations (\ref{eq:SO3_eoms_example}) and the
boundary conditions (\ref{eq:SO3_bcs_example}). The Lagrange multipliers
$\Lambdab_1\in\so^*(3), ~\Lambdab_2\in\so(3)$ correspond to the kinematic and dynamics constraints (\ref{eq:SO3_eoms_example}), respectively.
\end{thm}

\begin{rem} Note that the equations of motion that
arise from the Lagrange--d'Alembert principle are used to define the
dynamic constraints. In effect, we minimize $\Jcal$
subject to satisfying the Lagrange--d'Alembert principle. Analogously, the discrete version of the
Lagrange--d'Alembert principle will be used to derive the discrete
equations of motion in the discrete optimal control problem to be
studied in Section \ref{subsec:Discrete_Variational}. This view is
in line with the approach in \cite{Junge:05} in that we do not
discretize the equations of motion directly, but, instead, we
discretize the Lagrange--d'Alembert principle. These two approaches
are not equivalent in general.\end{rem}

\begin{cor} The necessary optimality conditions of Theorem \ref{thm:SO3_example2} are equivalent to those of
Theorem \ref{thm:SO3_example}.\end{cor}

\begin{pf} In Theorem \ref{thm:SO3_example}, differentiate $\Lambdab_2$ once and then use all three
differential equations to replace $\Lambdab_1$ and $\Lambdab_2$ with
expressions involving only $\taub, ~\Mb$ and $\Omegab$.\end{pf}

\section{Discrete-Time Results}\label{sec:discrete_time}

\subsection{Problem Formulation}\label{subsec:formulation_discrete_example}

In this section we give the discrete version of the problem
introduced in Section \ref{subsec:formulation_continuous_example}.
So, we consider minimizing the norm squared of the control torque
$\taub$ subject to satisfaction of the discrete Lagrange--d'Alembert
principle for the rigid body whose configuration and body angular
velocity at time step $t_k$ are given by $\Rb_k\in\SO$ and
$\Omegab_k\in\so(3)$, respectively. The kinematic constraint may be
expressed as \be\label{eq:discrete_RB_kinematics}
  \Rb_{k+1} = \Rb_k\exp\left(h\Omegab_k\right)=\Rb_k\gb_k,
\ee where $h$ is the integration time step,
$\exp:\so(3)\rightarrow\SO$ is the exponential map and
$\gb_k=\exp(h\Omegab_k)$. The boundary conditions are given by
$(\Rb^*_0,\Omegab^*_0)$ and $(\Rb^*_N,\Omegab^*_{N-1})$, where
$t_0=0$ and $N=T/h$ is such that $t_N=T$.

The reason we constrain $\Omegab$ at $t=h(N-1)$ instead of at $t=hN$
is that a constraint on $\Omegab_k\in\so(3)$ corresponds, by left
translations to a constraint on $\dot{\Rb}_k\in\Tsf_{\Rb_k}\SO$. In
turn, in the discrete setting and depending on the choice of
discretization, this corresponds to a constraint on the neighboring
discrete points. With our
choice of discretization (equation
(\ref{eq:discrete_RB_kinematics})), this corresponds to constraints
on $\Rb_k$ and $\Rb_{k+1}$. Hence, to ensure that the effect of the
terminal constraint on $\Omegab$ is correctly accounted for, the
constraint must be imposed on $\Omegab_{N-1}$, which entails some
constraints on variations at both $\Rb_{N-1}$ and
$\Rb_N$. We will return to this point later in the paper.

Equation (\ref{eq:discrete_RB_kinematics}) is just one way of
discretizing the kinematics of the rigid body.
In the case of planar rigid body dynamics, this leads to the first-order Euler approximation. However, on $\SO$, our approach yields a novel discretization.
We make the above choice for discretization as it
guarantees, in general, that the angular velocity matrix $\Omegab_k$
remains on the algebra $\so(3)$ by using the exponential map. This
is natural to do in the context of discrete variational numerical
solvers (for both initial value and two point boundary value
problems). Following the methodology of \cite{Junge:05}, we have the following
optimal control problem.

\begin{prb}\label{pb:SO3_discrete_example0} Minimize
\be\label{eq:SO3_cost_LD_discrete_example0}
  \Jcal &=&
  \sum_{k=0}^N\frac{1}{2}\ll\taub_k,\taub_k\gg_*
\ee subject to
\begin{enumerate}
    \item satisfying the \emph{discrete Lagrange--d'Alembert principle}:
\be\label{eq:SO3_LD_discrete_example0}
      \deltab\sum_{k=0}^{N-1}\frac{1}{2}\left<\Jb\left(\Omegab_k\right),\Omegab_k\right>+\sum_{k=0}^N\left<\taub_k,\Wb_k\right>=0,
    \ee
    subject to $\Rb_0=\Rb_0^*$, $\Rb_N=\Rb_N^*$ and $\Rb_{k+1} =\Rb_k\gb_k$,
    $k=0,1,\ldots,N-1$, where $\Wb_k$ is the variation vector field at time step $t_k$ satisfying
    $\deltab\Rb_k=\Rb_k\Wb_k$,
    \item and the \emph{boundary conditions} \be\label{eq:SO3_bcs_discrete_example0}
    \Rb_0&=&\Rb_0^*, ~\Omegab_0=\Omegab^*_0, \nn\\
    \Rb_N&=&\Rb_N^*, ~\Omegab_{N-1}=\Omegab^*_{N-1}.
\ee
\end{enumerate}
\end{prb}

The following formulation is equivalent, where the discrete rigid body equations of
motion replace the Lagrange--d'Alembert principle constraint.

\begin{prb}\label{pb:SO3_discrete_example} Minimize
\be\label{eq:SO3_cost_LD_discrete_example}
  \Jcal &=&
  \sum_{k=0}^N\frac{1}{2}\ll\taub_k,\taub_k\gg_*
\ee subject to
\begin{enumerate}
    \item the \emph{discrete dynamics} \begin{align}\label{eq:SO3_eoms_discrete_example}
      &\Rb_{k+1} = \Rb_k\gb_k, ~k=0,\ldots,N-1\nn\\
      &\Mb_k =
      \Ad_{\gb_k}^*\left(h\taub_k+\Mb_{k-1}\right), ~k=1,\ldots,N-1,\nn\\
      &\Mb_k=\Jb\left(\Omegab_k\right), ~k=0,\ldots,N-1,
    \end{align}
    \item and the \emph{boundary conditions} \be\label{eq:SO3_bcs_discrete_example}
    \Rb_0&=&\Rb^*_0, ~\Omegab_0=\Omegab^*_0, \nn\\
    \Rb_N&=&\Rb^*_N,
    ~\Omegab_{N-1}=\Omegab^*_{N-1}.
\ee
\end{enumerate}
\end{prb}

Regarding terminal velocity conditions, note that
in the second of equations (\ref{eq:SO3_eoms_discrete_example}) if
we let $k=N$ we find that $\Omegab_N$ appears in the equation. A constraint on $\Omegab_N$ dictates
constraints at the points $\Rb_N$ and $\Rb_{N+1}$ through the first
equation in (\ref{eq:SO3_eoms_discrete_example}). Since we only
consider time points up to $t=Nh$, we can not allow $k=N$ in the
second of equations (\ref{eq:SO3_eoms_discrete_example}) and hence
our terminal velocity constraints are posed in terms of
$\Omegab_{N-1}$ instead of $\Omegab_N$.

As mentioned above, $\Wb_k$ is a variation vector field associated
with the perturbed group element $\Rb_k^{\epsilon}$. Likewise, we
need to define a variation vector field associated with the element
$\gb_k=\exp(h\Omegab_k)$. First, let the perturbed variable
$\gb_k^{\epsilon}$ be defined by \be\label{eq:variations_in_gk}
\gb_k^{\epsilon} = \gb_k \exp(\epsilon h \deltab \Omegab_k),\ee
where \be\deltab\Omegab=\frac{\partial\Omegab_k^{\epsilon}}{\partial
\epsilon}\bigg|_{\epsilon=0}.\nn\ee Note that
$\gb_k^{\epsilon}\big|_{\epsilon=0}=\gb_k$ as desired. Moreover, we
have \be\label{eq:Deltab_g} \deltab \gb_k 
  = \gb_k (h\deltab\Omega_k) \exp(\epsilon h \delta \Omega_k)\big|_{\epsilon=0}= h\gb_k
  \deltab\Omegab_k.\ee This will be needed later when taking
  variations.

\subsection{The Discrete Lagrange--d'Alembert Principle and the Rigid
Body Equations of Motion}\label{sec:discrete_lag_RB}

In this section we derive the discrete forced rigid body equations
of motion (equations (\ref{eq:SO3_eoms_discrete_example})) starting
with the discrete Lagrange--d'Alembert principle. As in the continuous case, we will compute the constrained variation of $\deltab\Omegab_k$. We begin by
rewriting the kinematic constraint as
$\exp^{-1}\left(\Rb^{-1}_k\Rb_{k+1}\right)=h\Omegab_k$, which
is easier to handle as an expression over the Lie algebra. Take variations to obtain, $-\Rb_k^{-1}\left(\deltab\Rb_k\right)\Rb^{-1}_k\Rb_{k+1}+\Rb^{-1}_k\deltab\Rb_{k+1}=h\gb_k\cdot\deltab\Omegab_k$, which is equivalent to
$-\Wb_k\gb_k+\gb_k\Wb_{k+1}=h\gb_k\deltab\Omegab_k$, or
\be\label{eq:SO3_variations_example_discrete}
\deltab\Omegab_k=\frac{1}{h}\left[-\Ad_{\gb_k^{-1}}\Wb_k+\Wb_{k+1}\right].
\ee Note that this is an expression over the Lie algebra $\so(3)$.

Taking direct variations of the
cost functional we obtain \begin{align}
  &0=
  \left<\taub_0-\frac{1}{h}\Ad_{\gb_0^{-1}}^*\Jb\left(\Omegab_{0}\right),\Wb_0\right>\nn\\
  &+\left<\taub_N+\frac{1}{h}\Jb\left(\Omega_{N-1}\right),\Wb_N\right>\nn\\
  &+\sum_{k=1}^{N-1}\left<\taub_k
  -\frac{1}{h}\Ad_{\gb_k^{-1}}^*\Jb\left(\Omegab_{k}\right)+\frac{1}{h}\Jb\left(\Omegab_{k-1}\right),\Wb_k\right>.\nn
\end{align} where we have used equation
(\ref{eq:SO3_variations_example_discrete}). 
By the boundary conditions $\Rb_0=\Rb_0^*$ and $\Rb_N=\Rb_N^*$,
we have $\Wb_0=0$ and
$\Wb_N=0$. Since $\Wb_k$, $k=1,\ldots,N-1$,
are arbitrary and independent, we obtain the equivalent of equations
(\ref{eq:SO3_eoms_discrete_example}).

\subsection{Discrete-Time Variational Optimal Control Problem}\label{subsec:Discrete_Variational}

Analogous to the direct
approach in continuous time, here we derive the necessary optimality
conditions in a form that does not involve the use of Lagrange
multipliers. Using equation
(\ref{eq:SO3_variations_example_discrete}) and taking variation of
the second of equations (\ref{eq:SO3_eoms_discrete_example}), we
obtain
\be\label{eq:variation_in_Tau}\deltab\taub_k&=&\Ad^*_{\gb_k^{-1}}\bigg(\frac{1}{h^2}\Jb\left(\Wb_{k+1}-\Ad_{\gb_k^{-1}}\Wb_k\right)\nn
\\&&
+\frac{1}{h}\left[\Wb_{k+1}-\Ad_{\gb_k^{-1}}\Wb_k,\Jb\left(\Omegab_k\right)\right]\bigg)\nn\\&&
-\frac{1}{h^2}\Jb\left(\Wb_k-\Ad_{\gb_{k-1}^{-1}}\Wb_{k-1}\right),
\ee for $k=1,\ldots,N-1$. Taking variations of the cost
functional (\ref{eq:SO3_cost_LD_discrete_example}) and using equation (\ref{eq:variation_in_Tau}), one obtains after a
tedious but straightforward computation an expression for
$\deltab\Jcal$ in terms of $\deltab\taub_k$, which we omit because of space restrictions. 
When $\deltab\Jcal$ is equated to zero, the resulting equation gives (boundary) conditions on
$\taub_0,\taub_1,\taub_{N-1},\taub_{N}$ as well as discrete
evolution equations that are written in algebraic nonlinear form
as:\be\label{eq:DiscreteDirect}0&=&-\frac{1}{h^2}\bigg(\Jb\left(\taub_k^{\sharp}\right)-\Ad^*_{\gb_k^{-1}}\Jb\left(\taub_{k+1}^{\sharp}\right)
\nn\\&&-\Jb\left(\Ad_{\gb_{k-1}^{-1}}\taub_{k-1}^{\sharp}\right)
+\Ad^*_{\gb_{k}^{-1}}\Jb\left(\Ad_{\gb_{k}^{-1}}\taub_{k}^{\sharp}\right)\bigg)
\nn\\&&-\frac{1}{h}\bigg(\Ad^*_{\gb_{k}^{-1}}\left[\Jb\left(\Omegab_{k}\right),\Ad_{\gb_{k}^{-1}}\left(\taub_{k}^{\sharp}\right)\right]
\nn\\&&-\frac{1}{h}\left[\Jb\left(\Omegab_{k-1}\right),\Ad_{\gb_{k-1}^{-1}}\left(\taub_{k-1}^{\sharp}\right)\right]\bigg),\ee
for $k=2,\ldots,N-2$.

A Lagrange multiplier approach yields the following equivalent theorem.

\begin{thm}\label{thm:lag_disc_so3} The necessary optimality conditions for the discrete
Problem \ref{pb:SO3_discrete_example} are \be\label{eq:lag_disc_so3}
\Rb_{k+1} &=& \Rb_k\gb_k, ~k=1,\ldots,N-2\nn\\
      \Mb_k &=& \Ad_{\gb_k}^*\left(h\taub_k+\Mb_{k-1}\right),
      ~k=1,\ldots,N-1\nn\\
  0 &=& \Lambdab_{k-1}^1-\Ad_{\gb_k^{-1}}^*\Lambdab_k^1, ~k=2,\ldots,N-2\nn\\
  0
  &=&-\Lambdab_k^1+\Jb\left(\Lambdab_k^2\right)-\Jb\left(\Ad_{\gb_{k+1}}\Lambdab_{k+1}^2\right)\\&&+
  h\left[\Mb_k,\Lambdab_k^2\right], ~k=1,\ldots,N-2\nn\\
\taub_k &=& h\left(\Ad_{\gb_k}\Lambdab_k^2\right)^{\flat}, ~k=1,\ldots,N-1\nn  \\
  \Mb_k&=&\Jb\left(\Omegab_k\right), ~k=0,\ldots,N-1,\nn
\ee and the boundary conditions \be
  \Rb_0&=&\Rb_0^*, ~\Rb_1=\Rb_0^*\gb_0^*, ~\Omegab_0=\Omegab_0^*\nn\\
  \Rb_N&=&\Rb_N^*, ~\Rb_{N-1}=\Rb_N^*\left(\gb_{N-1}^*\right)^{-1}, ~\Omegab_{N-1}=\Omegab_{N-1}^*\nn\\
  \taub_0&=&\taub_N=0,\nn
\ee where $\gb_0^*=\exp(h\Omegab_0^*)$ and
$\gb_{N-1}^*=\exp\left(h\Omegab_{N-1}^*\right)$.\end{thm}

\section{Numerical Approach and Results}\la{sec:numerics}

The first-order optimality equations,
equation (\ref{eq:DiscreteDirect}), in combination with the boundary
conditions,
\[ \Rb_0= \Rb_0^*,\;\ \Rb_N= \Rb_N^*,\;\ \Omegab_0=\Omegab_0^*,\;\ \mbox{and  }
\Omegab_{N-1}=\Omegab_{N-1}^*, \] leave the torques
$\taub_1,\ldots,\taub_{N-1}$, and the angular velocities
$\Omegab_1,\ldots,\Omegab_{N-2}$ as unknowns. By substituting the
relations $\gb_k=\exp(h\Omegab_k)$,
$\mathbf{M}_k=\mathbf{J}(\Omegab_k)$, we can rewrite the necessary
conditions (\ref{eq:DiscreteDirect}) as follows,

\begin{align*}
0=&-\frac{1}{h^2}\Big(\Jb(\taub^\sharp_k)-\Ad^*_{\exp(-h\Omegab_k)}\Jb(\taub^\sharp_{k+1})\\
&\qquad\qquad-\Jb(\Ad_{\exp(-h\Omegab_{k-1})}\taub^\sharp_{k-1}\\
&\qquad\qquad+\Ad^*_{\exp(-h\Omegab_k)}\Jb(\Ad_{\exp(-h\Omegab_k)}\taub^\sharp_k)\Big)\\
&-\frac{1}{h}\Big(\Ad^*_{\exp(-h\Omegab_k)}\Big[\Jb(\Omegab_k),\Ad_{\exp(-h\Omegab_k)}(\taub^\sharp_k)\Big]\\
&\qquad\qquad-\frac{1}{h}\Big[\Jb(\Omegab_{k-1}),\Ad_{\exp(-h\Omegab_{k-1})}(\taub^\sharp_{k-1})\Big]\Big),\\
\intertext{where $k=2,\ldots,N-2$, and the discrete evolution equations, given by line 2 of (\ref{eq:lag_disc_so3}), can be written as}
0=&\Jb(\Omegab_k)-\Ad^*_{\exp(h\Omegab_k)} (h\taub_k+\mathbf{J}(\Omegab_{k-1})),\\
\intertext{where $k=1,\ldots,N-1$. In addition, we use the boundary
conditions on $\Rb_0$ and $\Rb_N$, together with the update step given by line 1 of (\ref{eq:lag_disc_so3}) to give the last constraint,}
0=&\log\Big(\Rb_N^{-1} \Rb_0 \exp(h\Omegab_0)\ldots
\exp(h\Omegab_{N-1})\Big),
\end{align*}
where $\log$ is the logarithm map on $\SO$.

At this point it should be noted that one important advantage of the
manner in which we have discretized the optimal control problem is
that it is $\SO$-equivariant. This is to say that if we rotated all
the boundary conditions by a fixed rotation matrix, and solved the
resulting discrete optimal control problem, the solution we would
obtain would simply be the rotation of the solution of the original
problem. This can be seen quite clearly from the fact that the
discrete problem is expressed in terms of body coordinates, both in
terms of body angular velocities and body forces. In addition, the
initial and final attitudes $\Rb_0$ and $\Rb_N$ only enter in the
last equation as a relative rotation.

The $\SO$-equivariance of our numerical method is desirable, since
it ensures that our results are independent of the choice of
coordinate frames. This is in contrast to methods based on
coordinatizing the rotation group using quaternions and Euler
angles.

The equations above take values in $\mathfrak{so}(3)$.
Consider the Lie algebra isomorphism between $\mathbb{R}^3$ and
$\mathfrak{so}(3)$ given by,
\[\mathbf{v}=(v_1,v_2,v_3)\mapsto\hat{\mathbf{v}}=\begin{bmatrix}0 & -v_3 & v_2\\ v_3 & 0 & -v_1\\-v_2 & v_1 & 0\end{bmatrix},\]
which maps 3-vectors to $3\times 3$ skew-symmetric matrices.
In particular, we have the following identities,
\begin{align*}
[\hat{\mathbf{u}},\hat{\mathbf{v}}]=(\mathbf{u}\times\mathbf{v})\,\hat{}\,,\quad
\Ad_\mathbf{A}\hat{\mathbf{v}}=(\mathbf{A}\mathbf{v})\,\hat{}\,.
\end{align*}
Furthermore, we identify $\mathfrak{so}(3)^*$ with $\mathbb{R}^3$ by
the usual dot product, that is to say if $\mathbf{\Pi}$,
$\mathbf{v}\in\mathbb{R}^3$, then $\langle
\mathbf{\Pi},\hat{\mathbf{v}}\rangle=\mathbf{\Pi}\cdot\mathbf{v}$.
With this identification, we have that
$\Ad^*_{\mathbf{A}^{-1}}\mathbf{\Pi}=\mathbf{A}\mathbf{\Pi}$. Using
the identities above, we write the necessary conditions using
matrix-vector products and cross products. Then, each of the
equations can be interpreted as $3$-vector valued functions, and the
system of equations can be considered as a $3(2N-3)$-vector valued
function, which is precisely the dimensionality of the unknowns.
This reduces the discrete optimal problem to a nonlinear root
finding problem.

We used a Newton-Armijo method, a line search algorithm
that uses the Newton search direction, and backtracking to
ensure sufficient descent of the residual error. The Jacobian is
constructed column by column, where the $k$-th column is computed
using the following approximation,
\[ \frac{\partial\mathbf{F}}{\partial x_k}(\mathbf{x}) = \frac{1}{\epsilon}\operatorname{Im}[{\mathbf{F}(\mathbf{x}+i\epsilon\mathbf{e}_k})],\]
where $i=\sqrt{-1}$, $\mathbf{e}_k$ is a basis vector in the $x_k$
direction, and $\epsilon$ is of the order of machine epsilon. This
method is more accurate than a finite-difference approximation as it
does not suffer from round-off errors.

\begin{figure}[h!]
\begin{center}
\includegraphics[width=3in, bb=20 20 500 380, clip=true]{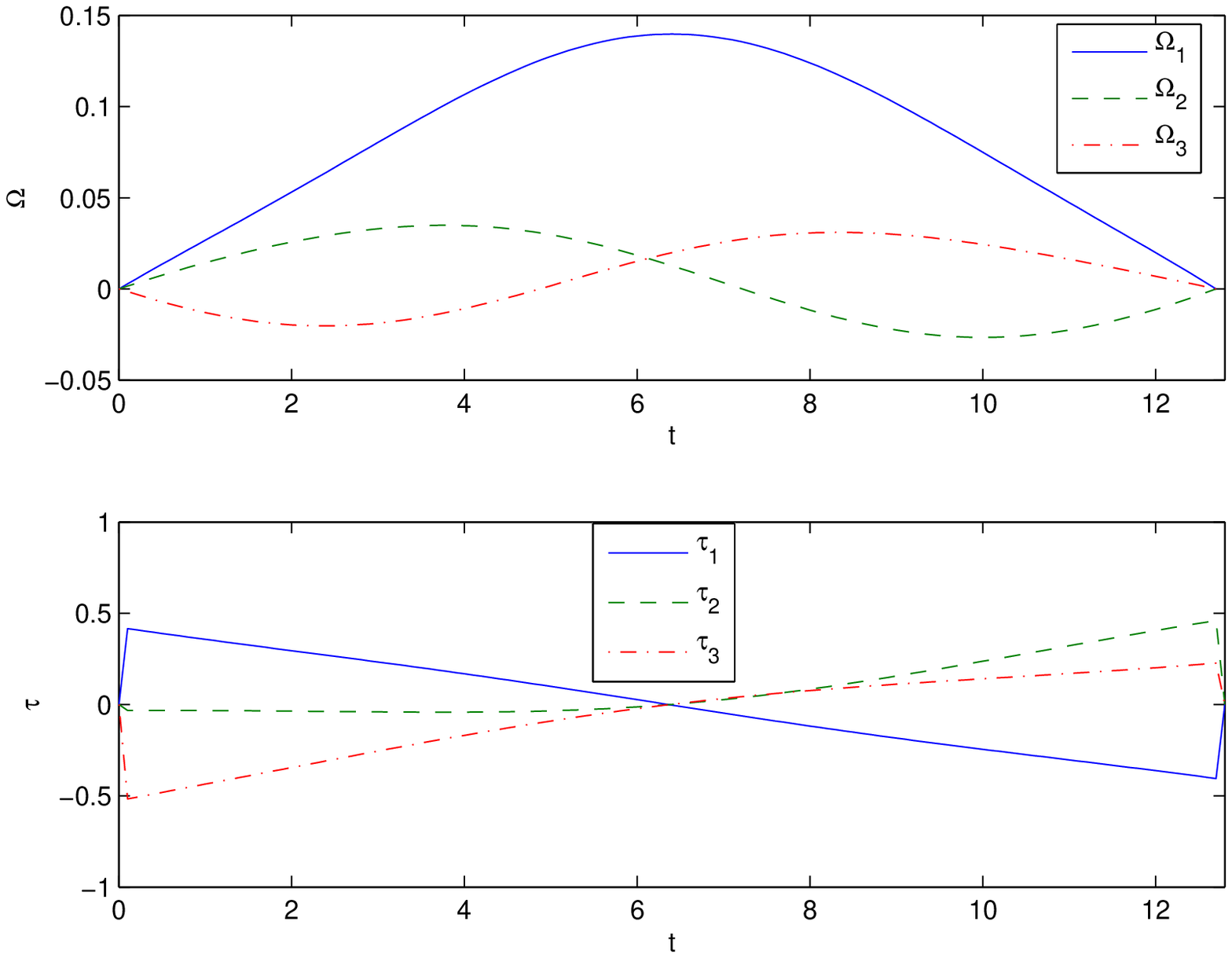}
\vspace{-5pt}
\caption{Discrete optimal rest-to-rest maneuver in
$\SO$.}\label{fig:result_so3}
\vspace*{-10pt}
\end{center}
\end{figure}

In our numerical simulation, we computed an optimal trajectory for a
rest-to-rest maneuver, as illustrated in Figure (\ref{fig:result_so3}). In this
simulation $N=128$, and essentially identical results were obtained
for $N=64$. It is worth noting that the results are not rotationally
symmetric about the midpoint of the simulation interval, which is
due to the fact that our choice of update,
$\Rb_{k+1}=\Rb_k\exp(h\Omegab_k)$, does not exhibit time-reversal
symmetry. In a forthcoming publication, we will introduce a
reversible algorithm to address this issue.

\begin{figure}[h!]
\begin{center}
\includegraphics[width=3in, bb=20 20 500 380, clip=true]{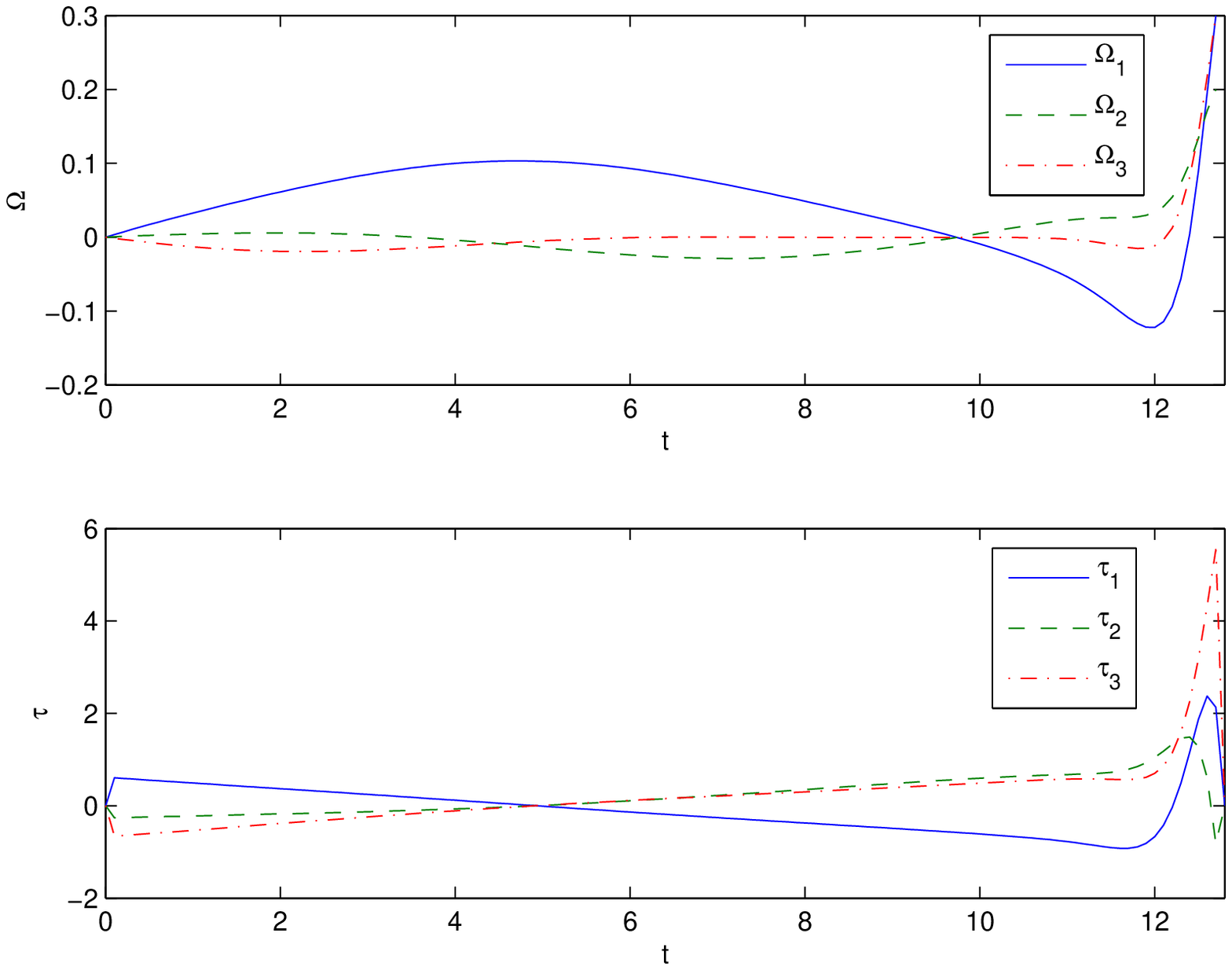}
\vspace*{-5pt}
\caption{Discrete optimal slew-up maneuver in
$\SO$.}\label{fig:slewup_so3}
\vspace*{-10pt}
\end{center}
\end{figure}

We also present results for an optimal slew-up maneuver in $T=12.8$ seconds from 
zero initial angular velocity to a final angular velocity of 
$\Omega_T=[0.3\;0.2\;0.3]^{\Trm}$, illustrated in Figure (\ref{fig:slewup_so3}). The 
resolution is $N=128$, and essentially identical results were obtained for $N=64$.

\section{Conclusion}\la{sec:conclusion}

In this paper we studied the continuous- and discrete-time optimal
control problem for the rigid body, where the cost to be minimized
is the external torque applied to move the rigid body from an
initial condition to some pre-specified terminal condition. In the
discrete setting, we use the discrete Lagrange--d'Alembert principle
to obtain the discrete equations of motion. The kinematics were
discretized to guarantee that the flow in phase space remains on the
Lie group $\SO$ and its algebra $\so(3)$. We described how the
necessary conditions can be solved for the general three-dimensional
case and gave a numerical example for a three-dimensional rigid body
maneuver.

Currently, we are investigating the connections with Pontryagin's maximum principle in continuous- and discrete-time.
Additionally, we wish to generalize the result to general Lie groups with applications other than rigid body motion on $\SO$.
In particular, we are interested in controlling the motion of a rigid body \emph{in space},
which corresponds to motion on the non-compact Lie group $\SEthree$.

\section*{Acknowledgments}
The research of Melvin Leok was partially supported by NSF grant DMS-0504747 and a University of Michigan Rackham faculty grant.
The research of Anthony Bloch was supported by NSF grants DMS-0305837, CMS-0408542, and DMS-0604307.

\bibliographystyle{IEEEtran}
\bibliography{bibliography}

\end{document}